\def\Resetstrings{%Clears all strings before processing reference listing.
    \def\present{ }\let\bgroup={\let\egroup=}%primitive TeX
    \def\Astr{}\def\astr{}\def\Atest{}\def\atest{}%
    \def\Bstr{}\def\bstr{}\def\Btest{}\def\btest{}%
    \def\Cstr{}\def\cstr{}\def\Ctest{}\def\ctest{}%
    \def\Dstr{}\def\dstr{}\def\Dtest{}\def\dtest{}%
    \def\Estr{}\def\estr{}\def\Etest{}\def\etest{}%
    \def\Fstr{}\def\fstr{}\def\Ftest{}\def\ftest{}%
    \def\Gstr{}\def\gstr{}\def\Gtest{}\def\gtest{}%
    \def\Hstr{}\def\hstr{}\def\Htest{}\def\htest{}%
    \def\Istr{}\def\istr{}\def\Itest{}\def\itest{}%
    \def\Jstr{}\def\jstr{}\def\Jtest{}\def\jtest{}%
    \def\Kstr{}\def\kstr{}\def\Ktest{}\def\ktest{}%
    \def\Lstr{}\def\lstr{}\def\Ltest{}\def\ltest{}%
    \def\Mstr{}\def\mstr{}\def\Mtest{}\def\mtest{}%
    \def\Nstr{}\def\nstr{}\def\Ntest{}\def\ntest{}%
    \def\Ostr{}\def\ostr{}\def\Otest{}\def\otest{}%
    \def\Pstr{}\def\pstr{}\def\Ptest{}\def\ptest{}%
    \def\Qstr{}\def\qstr{}\def\Qtest{}\def\qtest{}%
    \def\Rstr{}\def\rstr{}\def\Rtest{}\def\rtest{}%
    \def\Sstr{}\def\sstr{}\def\Stest{}\def\stest{}%
    \def\Tstr{}\def\tstr{}\def\Ttest{}\def\ttest{}%
    \def\Ustr{}\def\ustr{}\def\Utest{}\def\utest{}%
    \def\Vstr{}\def\vstr{}\def\Vtest{}\def\vtest{}%
    \def\Wstr{}\def\wstr{}\def\Wtest{}\def\wtest{}%
    \def\Xstr{}\def\xstr{}\def\Xtest{}\def\xtest{}%
    \def\Ystr{}\def\ystr{}\def\Ytest{}\def\ytest{}%
}
\Resetstrings

\def\Refformat{%Determines the kind of reference by the presence or
         \if\Jtest\present
             {\if\Vtest\present\journalarticleformat
                  \else\conferencereportformat\fi}
            \else\if\Btest\present\bookarticleformat
               \else\if\Rtest\present\technicalreportformat
                  \else\if\Itest\present\bookformat
                     \else\otherformat\fi\fi\fi\fi}

\def\Rpunct{%Default punctuation control strings if the punctuation
   \def\Lspace{ }%
   \def\Lperiod{ }%  .
   \def\Lcomma{ }%    ,
   \def\Lquest{ }%     ?
   \def\Lcolon{ }%   :
   \def\Lscolon{ }%   ;
   \def\Lbang{ }%      !
   \def\Lquote{ }%   '
   \def\Lqquote{ }%   "
   \def\Lrquote{ }%    `
   \def\Rspace{}%
   \def\Rperiod{.}%  .
   \def\Rcomma{,}%    ,
   \def\Rquest{?}%     ?
   \def\Rcolon{:}%   :
   \def\Rscolon{;}%   ;
   \def\Rbang{!}%      !
   \def\Rquote{'}%   '
   \def\Rqquote{"}%   "
   \def\Rrquote{`}%    `
   }

\def\Lpunct{%Default punctuation control strings if the punctuation
   \def\Lspace{}%
   \def\Lperiod{\unskip.}%  .
   \def\Lcomma{\unskip,}%    ,
   \def\Lquest{\unskip?}%     ?
   \def\Lcolon{\unskip:}%   :
   \def\Lscolon{\unskip;}%   ;
   \def\Lbang{\unskip!}%      !
   \def\Lquote{\unskip'}%   '
   \def\Lqquote{\unskip"}%   "
   \def\Lrquote{\unskip`}%    `
   \def\Rspace{\spacefactor=1000}%
   \def\Rperiod{\spacefactor=3000}%  .
   \def\Rcomma{\spacefactor=1250}%    ,
   \def\Rquest{\spacefactor=3000}%     ?
   \def\Rcolon{\spacefactor=2000}%   :
   \def\Rscolon{\spacefactor=1250}%   ;
   \def\Rbang{\spacefactor=3000}%      !
   \def\Rquote{\spacefactor=1000}%   '
   \def\Rqquote{\spacefactor=1000}%   "
   \def\Rrquote{\spacefactor=1000}%    `
   }

\def\Refstd{%Standard control strings for formatting bibliography listings.
     \def\Acomma{\unskip, }%between multiple author names
     \def\Aand{\unskip\ and }%between two author names
     \def\Aandd{\unskip\ and }%between last two of multiple author names
     \def\Ecomma{\unskip, }%between multiple editor names
     \def\Eand{\unskip\ and }%between two editor names
     \def\Eandd{\unskip\ and }%between last two of multiple author names
     \def\acomma{\unskip, }%same for authors of reviewed material
     \def\aand{\unskip\ and }%same for authors of reviewed material
     \def\aandd{\unskip\ and }%same for authors of reviewed material
     \def\ecomma{\unskip, }%same for translators
     \def\eand{\unskip\ and }%same for translators
     \def\eandd{\unskip\ and }%same for translators
     \def\Namecomma{\unskip, }%same for citations using authors' names
     \def\Nameand{\unskip\ and }%same for citations using authors' names
     \def\Nameandd{\unskip\ and }%same for citations using authors' names
     \def\Revcomma{\unskip, }%between last and first name of reversed name
     \def\Initper{.\ }%punctuation after initial
     \def\Initgap{\dimen0=\spaceskip\divide\dimen0 by 2\hskip-\dimen0}%
   }

\def\Smallcapsaand{%Smallcaps redefinition of \Aand and \Aandd for \Refstd
     \def\Aand{\unskip\bgroup{\Smallcapsfont\ AND }\egroup}%
     \def\Aandd{\unskip\bgroup{\Smallcapsfont\ AND }\egroup}%
     \def\eand{\unskip\bgroup\Smallcapsfont\ AND \egroup}%
     \def\eandd{\unskip\bgroup\Smallcapsfont\ AND \egroup}%
   }

\def\Smallcapseand{%Smallcaps redefinition of \Eand, \Eeand, etc for Refstd
     \def\Eand{\unskip\bgroup\Smallcapsfont\ AND \egroup}%
     \def\Eandd{\unskip\bgroup\Smallcapsfont\ AND \egroup}%
     \def\aand{\unskip\bgroup\Smallcapsfont\ AND \egroup}%
     \def\aandd{\unskip\bgroup\Smallcapsfont\ AND \egroup}%
   }

   \def\Citefont{}%citations
   \def\ACitefont{}%alternate citations
   \def\Authfont{}%authors
   \def\Titlefont{}%titles
   \def\Tomefont{\sl}%journals or books
   \def\Volfont{}%volume number of journal
   \def\Flagfont{}%citation flag
   \def\Reffont{\rm}%set at beginning of reference listing
   \def\Smallcapsfont{\sevenrm}%small caps
   \def\Flagstyle#1{\hangindent\parindent\indent\hbox to0pt%flag style
       {\hss[{\Flagfont#1}]\kern.5em}\ignorespaces}%        for references

\def\Citebrackets{\Rpunct%defaults for putting citations in brackets [].
   \def\Lcitemark{\def\Cfont{\Citefont}[\bgroup\Cfont}%mark at left of citation
   \def\Rcitemark{\egroup]}%mark at right of citation
   \def\LAcitemark{\def\Cfont{\ACitefont}\bgroup\ACitefont}%
   \def\RAcitemark{\egroup}%mark at right of alternate citation
   \def\LIcitemark{\egroup}%mark at left of insertion in citation
   \def\RIcitemark{\bgroup\Cfont}%mark at right of insertion in citation
   \def\Citehyphen{\egroup--\bgroup\Cfont}%separater for string of citations
   \def\Citecomma{\egroup,\hskip0pt\bgroup\Cfont}%
   \def\Citebreak{}%mark between parts of citation (e.g. author\Citebreak date)
   }

\def\Citeparen{\Rpunct%defaults for putting citations in parenthesis ().
   \def\Lcitemark{\def\Cfont{\Citefont}(\bgroup\Cfont}%mark at left of citation
   \def\Rcitemark{\egroup)}%mark at right of citation
   \def\LAcitemark{\def\Cfont{\ACitefont}\bgroup\ACitefont}%
   \def\RAcitemark{\egroup}%mark at right of alternate citation
   \def\LIcitemark{\egroup}%mark at left of insertion in citation
   \def\RIcitemark{\bgroup\Cfont}%mark at right of insertion in citation
   \def\Citehyphen{\egroup--\bgroup\Cfont}%separater for string of citations
   \def\Citecomma{\egroup,\hskip0pt\bgroup\Cfont}%
   \def\Citebreak{}%mark between parts of citation (e.g. author\Citebreak date)
   }

\def\Citesuper{\Lpunct%defaults for making superscript citations
   \def\Lcitemark{\def\Cfont{\Citefont}\raise1ex\hbox\bgroup\bgroup\Cfont}%
   \def\Rcitemark{\egroup\egroup}%mark at right of citation
   \def\LAcitemark{\def\Cfont{\ACitefont}\bgroup\ACitefont}%
   \def\RAcitemark{\egroup}%mark at right of alternate citation
   \def\LIcitemark{\egroup\egroup}%mark at left of insertion in citation
   \def\RIcitemark{\raise1ex\hbox\bgroup\bgroup\Cfont}%
   \def\Citehyphen{\egroup--\bgroup\Cfont}%separater for string of citations
   \def\Citecomma{\egroup,\hskip0pt\bgroup%
      \Cfont}%separater for multiple citations
   \def\Citebreak{}%mark between parts of citation (e.g. author\Citebreak date)
   }

\def\Citenamedate{\Rpunct%defaults for making name-date citations
   \def\Lcitemark{%mark at left of citation--also sets internal punctuation
      \def\Citebreak{\egroup\ [\bgroup\Citefont}%separater in citation
      \def\Citecomma{\egroup]; %between multiple citations
         \bgroup\let\uchyph=1\Citefont}(\bgroup\let\uchyph=1\Citefont}%
   \def\Rcitemark{\egroup])}%mark at right of citation
   \def\LAcitemark{%mark at left of alternate citation
      \def\Citebreak{\egroup\ [\bgroup\Citefont}\def\Citecomma{\egroup], %
         \bgroup\ACitefont }\bgroup\let\uchyph=1\ACitefont}%
   \def\RAcitemark{\egroup]}%mark at right of alternate citation
  \def\Citehyphen{\egroup--\bgroup\Citefont}%separater for string of citations
   \def\LIcitemark{\egroup}%mark at left of insertion in citation
   \def\RIcitemark{\bgroup\Citefont}%mark at right of insertion in citation
   }
\def\annotations{}
\def\Flagstyle#1{\par\noindent\hbox to 30pt{[#1]\hfil}%
\hangindent=30pt\hangafter=1}

\def\ifundefined#1{\expandafter\ifx\csname#1\endcsname\relax}
\ifundefined{thesis}\def\bye{\end{document}}\fi

\def\endauthor{:\kern6pt}

\Refstd\Citebrackets%set general formats for reference list and citations
\def\Citefont{\rm}\def\Titlefont{\Reffont}\def\Volfont{\bf}\def\Tomefont{\it}%redefine some fonts

\def\journalarticleformat{\Reffont\let\uchyph=1\parindent=1.25pc\def\Comma{}%

\sfcode`\.=1000\sfcode`\?=1000\sfcode`\!=1000\sfcode`\:=1000\sfcode`\;=1000\sfcode`\,=1000%\frenchspacing
                \par\vfil\penalty-200\vfilneg%\filbreak
      \if\Ftest\present\Flagstyle\Fstr\fi%

\if\Atest\present\bgroup\Authfont\Astr\egroup\def\Comma{\unskip\endauthor}\fi%
        \if\Ttest\present\Comma\bgroup\Titlefont\Tstr\egroup\def\Comma{, }\fi%
         \if\etest\present\hskip.2em(\bgroup\estr\egroup)\def\Comma{\unskip,
}\fi%
          \if\Jtest\present\Comma\bgroup\Tomefont\Jstr\/\egroup\def\Comma{,
}\fi%

\if\Vtest\present\if\Jtest\present\hskip.2em\else\Comma\fi\bgroup\Volfont\Vstr\egroup\def\Comma{, }\fi%
            \if\Dtest\present\hskip.2em(\bgroup\Dstr\egroup)\def\Comma{, }\fi%
             \if\Ptest\present\bgroup, \Pstr\egroup\def\Comma{, }\fi%
              \if\ttest\present\Comma\bgroup\Titlefont\tstr\egroup\def\Comma{,
}\fi%

\if\jtest\present\Comma\bgroup\Tomefont\jstr\/\egroup\def\Comma{, }\fi%

\if\vtest,\present\if\jtest\present\hskip.2em\else\Comma\fi\bgroup\Volfont\vstr\egroup\def\Comma{, }\fi%
                 \if\dtest\present\hskip.2em(\bgroup\dstr\egroup)\def\Comma{,
}\fi%
                  \if\ptest\present\bgroup, \pstr\egroup\def\Comma{, }\fi%
                   \if\Gtest\present{\Comma Gov't ordering no.
}\bgroup\Gstr\egroup\def\Comma{, }\fi%
                    \if\Mtest\present\Comma MR
\#\bgroup\Mstr\egroup\def\Comma{, }\fi%
                     \if\Otest\present{\Comma\bgroup\Ostr\egroup.}\else{.}\fi%
\annotations  %  MODIFICATION
\vskip5ptplus1ptminus1pt}%\smallskip (but 3 is now 5)  MODIFICATION

\def\conferencereportformat{\Reffont\let\uchyph=1\parindent=1.25pc\def\Comma{}%

\sfcode`\.=1000\sfcode`\?=1000\sfcode`\!=1000\sfcode`\:=1000\sfcode`\;=1000\sfcode`\,=1000%\frenchspacing
                \par\vfil\penalty-200\vfilneg%\filbreak
      \if\Ftest\present\Flagstyle\Fstr\fi%

\if\Atest\present\bgroup\Authfont\Astr\egroup\def\Comma{\unskip\endauthor}\fi%
        \if\Ttest\present\Comma\bgroup\Titlefont\Tstr\egroup\def\Comma{, }\fi%
         \if\Jtest\present\Comma\bgroup\Tomefont\Jstr\/\egroup\def\Comma{,
}\fi%
          \if\Ctest\present\Comma\bgroup\Cstr\egroup\def\Comma{, }\fi%
           \if\Dtest\present\hskip.2em(\bgroup\Dstr\egroup)\def\Comma{, }\fi%
            \if\Mtest\present\Comma MR \#\bgroup\Mstr\egroup\def\Comma{, }\fi%
             \if\Otest\present{\Comma\bgroup\Ostr\egroup.}\else{.}\fi%
\annotations  %  MODIFICATION
\vskip5ptplus1ptminus1pt}%\smallskip (but 3 is now 5)  MODIFICATION

\def\bookarticleformat{\Reffont\let\uchyph=1\parindent=1.25pc\def\Comma{}%

\sfcode`\.=1000\sfcode`\?=1000\sfcode`\!=1000\sfcode`\:=1000\sfcode`\;=1000\sfcode`\,=1000%\frenchspacing
                \par\vfil\penalty-200\vfilneg%\filbreak
      \if\Ftest\present\Flagstyle\Fstr\fi%

\if\Atest\present\bgroup\Authfont\Astr\egroup\def\Comma{\unskip\endauthor}\fi%
        \if\Ttest\present\Comma\bgroup\Titlefont\Tstr\egroup\def\Comma{, }\fi%
         \if\etest\present\hskip.2em(\bgroup\estr\egroup)\def\Comma{\unskip,
}\fi%
          \if\Btest\present\Comma in
\bgroup\Tomefont\Bstr\/\egroup\def\Comma{\unskip, }\fi%
           \if\otest\present\Comma\ \bgroup\ostr\egroup\def\Comma{, }\fi%
            \if\Etest\present\Comma\bgroup\Estr\egroup\unskip,
\ifnum\Ecnt>1eds.\else ed.\fi\def\Comma{, }\fi%
             \if\Stest\present\Comma\bgroup\Sstr\egroup\def\Comma{, }\fi%
              \if\Vtest\present\ \bgroup\bf\Vstr\egroup\def\Comma{, }\fi%
               \if\Ntest\present\Comma no. \bgroup\Nstr\egroup\def\Comma{,
}\fi%
                \if\Itest\present\Comma\bgroup\Istr\egroup\def\Comma{, }\fi%
                 \if\Ctest\present\ (\bgroup\Cstr\egroup)\def\Comma{, }\fi%
                  \if\Dtest\present\Comma\bgroup\Dstr\egroup\def\Comma{, }\fi%
                   \if\Ptest\present\Comma pp.\ \Pstr\def\Comma{, }\fi%
\if\ttest\present\Comma\bgroup\Titlefont\Tstr\egroup\def\Comma{, }\fi%
                     \if\btest\present\Comma in
\bgroup\Tomefont\bstr\egroup\def\Comma{, }\fi%
                       \if\atest\present\Comma\bgroup\astr\egroup\unskip,
\if\acnt\present eds.\else ed.\fi\def\Comma{, }\fi%
                        \if\stest\present\Comma\bgroup\sstr\egroup\def\Comma{,
}\fi%
                         \if\vtest\present\Comma vol.
\bgroup\vstr\egroup\def\Comma{, }\fi%
                          \if\ntest\present\Comma no.
\bgroup\nstr\egroup\def\Comma{, }\fi%

\if\itest\present\Comma\bgroup\istr\egroup\def\Comma{, }\fi%

\if\ctest\present\Comma\bgroup\cstr\egroup\def\Comma{, }\fi%

\if\dtest\present\Comma\bgroup\dstr\egroup\def\Comma{, }\fi%
                              \if\ptest\present\Comma\pstr\def\Comma{, }\fi%
                               \if\Gtest\present{\Comma Gov't ordering no.
}\bgroup\Gstr\egroup\def\Comma{, }\fi%
                                \if\Mtest\present\Comma MR
\#\bgroup\Mstr\egroup\def\Comma{, }\fi%

\if\Otest\present{\Comma\bgroup\Ostr\egroup.}\else{.}\fi%
\annotations  %  MODIFICATION
\vskip5ptplus1ptminus1pt}%\smallskip (but 3 is now 5)  MODIFICATION

\def\bookformat{\Reffont\let\uchyph=1\parindent=1.25pc\def\Comma{}%

\sfcode`\.=1000\sfcode`\?=1000\sfcode`\!=1000\sfcode`\:=1000\sfcode`\;=1000\sfcode`\,=1000%\frenchspacing
                \par\vfil\penalty-200\vfilneg%\filbreak
      \if\Ftest\present\Flagstyle\Fstr\fi%

\if\Atest\present\bgroup\Authfont\Astr\egroup\def\Comma{\unskip\endauthor}%

\else\if\Etest\present\bgroup\def\Eand{\Aand}\def\Eandd{\Aandd}\Authfont\Estr\egroup\unskip, \ifnum\Ecnt>1eds.\else ed.\fi\def\Comma{, }%

\else\if\Itest\present\bgroup\Authfont\Istr\egroup\def\Comma{, }\fi\fi\fi%

\if\Ttest\present\Comma\bgroup\Tomefont\Tstr\/\egroup\def\Comma{\unskip, }%

\else\if\Btest\present\Comma\bgroup\Titlefont\Bstr\/\egroup\def\Comma{\unskip,
}\fi\fi%
            \if\otest\present\Comma\ \bgroup\ostr\egroup\def\Comma{, }\fi%

\if\etest\present\hskip.2em(\bgroup\estr\egroup)\def\Comma{\unskip, }\fi%
              \if\Stest\present\Comma\bgroup\Sstr\egroup\def\Comma{, }\fi%
               \if\Vtest\present\ \bgroup\bf\Vstr\egroup\def\Comma{, }\fi%
                \if\Ntest\present\Comma no. \bgroup\Nstr\egroup\def\Comma{,
}\fi%
                 \if\Atest\present\if\Itest\present
                         \Comma\bgroup\Istr\egroup\def\Comma{\unskip, }\fi%
                      \else\if\Etest\present\if\Itest\present
                              \Comma\bgroup\Istr\egroup\def\Comma{\unskip,
}\fi\fi\fi%
                     \if\Ctest\present\ (\bgroup\Cstr\egroup)\def\Comma{, }\fi%
                      \if\Dtest\present\Comma\bgroup\Dstr\egroup\def\Comma{,
}\fi%
                      \if\Ptest\present\bgroup, pp.\ \Pstr\egroup\def\Comma{,
}\fi% MODIFICATION

\if\ttest\present\Comma\bgroup\Titlefont\tstr\egroup\def\Comma{, }%

\else\if\btest\present\Comma\bgroup\Titlefont\bstr\egroup\def\Comma{, }\fi\fi%

\if\stest\present\Comma\bgroup\sstr\egroup\def\Comma{, }\fi%
                           \if\vtest\present\Comma vol.
\bgroup\vstr\egroup\def\Comma{, }\fi%
                            \if\ntest\present\Comma no.
\bgroup\nstr\egroup\def\Comma{, }\fi%

\if\itest\present\Comma\bgroup\istr\egroup\def\Comma{, }\fi%

\if\ctest\present\Comma\bgroup\cstr\egroup\def\Comma{, }\fi%

\if\dtest\present\Comma\bgroup\dstr\egroup\def\Comma{, }\fi%
                                \if\Gtest\present{\Comma Gov't ordering no.
}\bgroup\Gstr\egroup\def\Comma{, }\fi%
                                 \if\Mtest\present\Comma MR
\#\bgroup\Mstr\egroup\def\Comma{, }\fi%

\if\Otest\present{\Comma\bgroup\Ostr\egroup.}\else{.}\fi%
\annotations  %  MODIFICATION
\vskip5ptplus1ptminus1pt}%\smallskip (but 3 is now 5)  MODIFICATION

\def\technicalreportformat{\Reffont\let\uchyph=1\parindent=1.25pc\def\Comma{}%

\sfcode`\.=1000\sfcode`\?=1000\sfcode`\!=1000\sfcode`\:=1000\sfcode`\;=1000\sfcode`\,=1000%\frenchspacing
                \par\vfil\penalty-200\vfilneg%\filbreak
      \if\Ftest\present\Flagstyle\Fstr\fi%

\if\Atest\present\bgroup\Authfont\Astr\egroup\def\Comma{\unskip\endauthor}%

\else\if\Etest\present\bgroup\def\Eand{\Aand}\def\Eandd{\Aandd}\Authfont\Estr\egroup\unskip, \ifnum\Ecnt>1eds.\else ed.\fi\def\Comma{, }%

\else\if\Itest\present\bgroup\Authfont\Istr\egroup\def\Comma{, }\fi\fi\fi%
          \if\Ttest\present\Comma\bgroup\Titlefont\Tstr\egroup\def\Comma{,
}\fi%
           \if\Atest\present\if\Itest\present
                   \Comma\bgroup\Istr\egroup\def\Comma{, }\fi%
                \else\if\Etest\present\if\Itest\present
                        \Comma\bgroup\Istr\egroup\def\Comma{, }\fi\fi\fi%
            \if\Rtest\present\Comma\bgroup\Rstr\egroup\def\Comma{, }\fi%
             \if\Ctest\present\ (\bgroup\Cstr\egroup)\def\Comma{, }\fi%
              \if\Dtest\present\Comma\bgroup\Dstr\egroup\def\Comma{, }\fi%
              \if\Ptest\present\bgroup, \Pstr\egroup\def\Comma{, }\fi%
               \if\ttest\present\Comma\bgroup\Titlefont\tstr\egroup\def\Comma{,
}\fi%
                \if\itest\present\Comma\bgroup\istr\egroup\def\Comma{, }\fi%
                 \if\rtest\present\Comma\bgroup\rstr\egroup\def\Comma{, }\fi%
                  \if\ctest\present\Comma\bgroup\cstr\egroup\def\Comma{, }\fi%
                   \if\dtest\present\Comma\bgroup\dstr\egroup\def\Comma{, }\fi%
                    \if\Gtest\present{\Comma Gov't ordering no.
}\bgroup\Gstr\egroup\def\Comma{, }\fi%
                     \if\Mtest\present\Comma MR
\#\bgroup\Mstr\egroup\def\Comma{, }\fi%
                      \if\Otest\present{\Comma\bgroup\Ostr\egroup.}\else{.}\fi%
\annotations  %  MODIFICATION
\vskip5ptplus1ptminus1pt}%\smallskip (but 3 is now 5)  MODIFICATION

\def\otherformat{\Reffont\let\uchyph=1\parindent=1.25pc\def\Comma{}%

\sfcode`\.=1000\sfcode`\?=1000\sfcode`\!=1000\sfcode`\:=1000\sfcode`\;=1000\sfcode`\,=1000%\frenchspacing
                \par\vfil\penalty-200\vfilneg%\filbreak
      \if\Ftest\present\Flagstyle\Fstr\fi%

\if\Atest\present\bgroup\Authfont\Astr\egroup\def\Comma{\unskip\endauthor}%

\else\if\Etest\present\bgroup\def\Eand{\Aand}\def\Eandd{\Aandd}\Authfont\Estr\egroup\unskip, \ifnum\Ecnt>1eds.\else ed.\fi\def\Comma{, }%

\else\if\Itest\present\bgroup\Authfont\Istr\egroup\def\Comma{, }\fi\fi\fi%
          \if\Ttest\present\Comma\bgroup\Titlefont\Tstr\egroup\def\Comma{,
}\fi%
            \if\Atest\present\if\Itest\present
                    \Comma\bgroup\Istr\egroup\def\Comma{, }\fi%
                 \else\if\Etest\present\if\Itest\present
                         \Comma\bgroup\Istr\egroup\def\Comma{, }\fi\fi\fi%
                 \if\Ctest\present\ (\bgroup\Cstr\egroup)\def\Comma{, }\fi%
                  \if\Dtest\present\Comma\bgroup\Dstr\egroup\def\Comma{, }\fi%
                  \if\Ptest\present\bgroup, \Pstr\egroup\def\Comma{, }\fi%
                   \if\Gtest\present{\Comma Gov't ordering no.
}\bgroup\Gstr\egroup\def\Comma{, }\fi%
                    \if\Mtest\present\Comma MR
\#\bgroup\Mstr\egroup\def\Comma{, }\fi%
                     \if\Otest\present{\Comma\bgroup\Ostr\egroup.}\else{.}\fi%
\annotations  %  MODIFICATION
\vskip5ptplus1ptminus1pt}%\smallskip (but 3 is now 5)  MODIFICATION
\catcode`\@=11
\def\matrixx#1{\null\,\vcenter{\normalbaselines\m@th
    \ialign{\hfil$\displaystyle ##$\hfil&&\quad\hfil$\displaystyle
##$\hfil\crcr    \mathstrut\crcr\noalign{\kern-\baselineskip}
    #1\crcr\mathstrut\crcr\noalign{\kern-\baselineskip}}}\,}
\catcode`\@=12

\def\annotations{}
\def\Flagstyle#1{\par\noindent\hbox to 30pt{[#1]\hfil}%
\hangindent=30pt\hangafter=1}
\def\widepost#1#2#3{\par\vspace{0.1in}\hbox{\vbox to #1bp{\vss%
     \special{" HPageCenter 0 #2 translate DEFAULTS #3 }}}\vspace{0.1in}}

\catcode`\@=11
\def\ATstartsection{\@startsection}
\def\resetATfont{\reset@font}
\def\zAT{\z@}
\def\RIfM@{\relax\protect\ifmmode}

\catcode`\@=12

\def\hexnumberZ#1{\ifnum#1<10 \number#1\else
 \ifnum#1=10 A\else\ifnum#1=11 B\else\ifnum#1=12 C\else
 \ifnum#1=13 D\else\ifnum#1=14 E\else\ifnum#1=15 F\fi\fi\fi\fi\fi\fi\fi}

\newcounter{alphactr}

\newcounter{romanctr}

\def\PLUS{+}
\def\TIMES{*}

\def\begindiagram{\def\normalbaselines{\baselineskip20pt
      \lineskip3pt \lineskiplimit3pt}}

\catcode`\@=11
\def\operatoratfont{\operator@font}
\let\@afterindentfalse\@afterindenttrue
\@afterindenttrue
\catcode`\@=12

\def\many#1{\manyxx#1 \cdots \manyxx#1}
\def\manyxx#1{\if#1\PLUS{+}\else\if#1\TIMES{*}\fi\fi}

\def\xhmode#1{\relax\ifmmode{#1}\else{\snap\hbox{$#1$}}\fi}
\def\snap{\hskip 0pt plus 4cm\penalty1000\hskip 0pt plus -4cm}
\def\PP{\xmode{\Bbb P\kern1pt}}
\def\email{jaffe{\kern0.5pt}@{\kern0.5pt}cpthree.unl.edu}
\def\C{{\Bbb C}\kern1pt}
\def\Q{{\Bbb Q}\kern1pt}
\def\F{{\Bbb F}\kern1pt}
\def\N{\xmode{\Bbb N}}
\def\Z{\xmode{\Bbb Z}}
\def\O{{\cal O}}
\def\Ker{\mathop{\operatoratfont Ker}\nolimits}
\def\Pic{\mathop{\operatoratfont Pic}\nolimits}
\def\Sing{\mathop{\operatoratfont Sing}\nolimits}
\def\mp[[#1||#2||#3]]{\snap\hbox{$#1:\kern3pt #2\ \to\ #3$}}
\def\mapx[[#1||#2]]{\snap\hbox{$#1\ \to\ #2$}}
\def\mapE#1{{\smash{\mathop{\longrightarrow}\limits^{#1}}}}
\def\smapE#1{{\smash{\mathop{\kern-10pt\rightarrow\kern-10pt}\limits^{#1}}}}
\def\diagramx#1{{$$\begindiagram\matrixx{#1}$$}}
\def\splitdiagram#1#2{
    \begin{flushleft}$\displaystyle\begindiagram\matrixx{#1}$\end{flushleft}
    \begin{flushright}$\displaystyle\begindiagram\matrixx{#2}$\end{flushright}}
\def\ses#1#2#3%
{  \setbox0=\hbox{$\matrixx{0&\mapE{}&#1&\mapE{}&#2&\mapE{}&#3&\mapE{}&0}$%
}  \ifdim\wd0<\hsize\diagramx{0&\mapE{}&#1&\mapE{}&#2&\mapE{}&#3&\mapE{}&
  0}\else\setbox0=\hbox{$\matrixx{0&\smapE{}&#1&\smapE{}&#2&\smapE{}&#3&
  \smapE{}&0}$}\ifdim\wd0<\hsize\diagramx{0&\smapE{}&#1&\smapE{}&#2&\smapE{}&
  #3&\smapE{}&0}\else\splitdiagram{0&\mapE{}&#1&\mapE{}&#2}{\mapE{}&#3
  &\mapE{}&0}\fi\fi}
\def\sesdot#1#2#3%
{  \setbox0=\hbox{$\matrixx{0&\mapE{}&#1&\mapE{}&#2&\mapE{}&#3&\mapE{}&0.}$%
}  \ifdim\wd0<\hsize\diagramx{0&\mapE{}&#1&\mapE{}&#2&\mapE{}&#3&\mapE{}&
  0.}\else\setbox0=\hbox{$\matrixx{0&\smapE{}&#1&\smapE{}&#2&\smapE{}&#3&
  \smapE{}&0.}$}\ifdim\wd0<\hsize\diagramx{0&\smapE{}&#1&\smapE{}&#2&\smapE{}&
  #3&\smapE{}&0.}\else\splitdiagram{0&\mapE{}&#1&\mapE{}&#2}{\mapE{}&#3
  &\mapE{}&0.}\fi\fi}
\def\o#1{\if#1\PLUS\oplus\else\if#1\TIMES\otimes\fi\fi}
\def\tL{{\tilde{L}}}
\def\tS{{\tilde{S}}}
\def\oM{{\overline{M}}}
\def\shC{{\cal{C}}}
\def\shL{{\cal{L}}}
\def\STCIs{set-theoretic complete intersections}
\def\WRT{with respect to}
\def\RHS{right hand side}
\def\IN{\subset}
\def\setof#1{\xmode{\{#1\}}}
\def\iso{\cong}
\def\qed{{\hfill$\square$}}
\def\red#1{#1_{\operatoratfont red}}
\def\half{{1\over2}}
\def\vec#1#2#3{#1_{#2}, \ldots, #1_{#3}}
\def\svec#1#2#3{#1_{#2} \many+ #1_{#3}}
\def\xmode#1{\relax\ifmmode{#1}\else{$#1$}\fi}
\def\pref#1{(\ref{#1})}
\def\block#1{\section{#1}}
\def\abs#1{{\vert{#1}\vert}}
\def\floor#1{\lfloor#1\rfloor}
\def\makeaddress{
     \vskip 0.15in
     \par\noindent {\footnotesize Department of Mathematics and Statistics,
                                  University of Nebraska}
     \par\noindent {\footnotesize Lincoln, NE 68588-0323, USA\ \ (\email)}}
\def\nsfatnebraska{ \def\thefootnote{\fnsymbol{footnote}}
     \par\noindent David B. Jaffe\protect\footnote{Partially supported by
                   the National Science Foundation.}
     \makeaddress\def\thefootnote{\arabic{footnote}}\setcounter{footnote}{0}}
\newenvironment{proof}{\trivlist \item[\hskip \labelsep{\sc
Proof.\kern1pt}]}{\endtrivlist}%\proof

\newenvironment{alphalist}{\begin{list}{(\alph{alphactr})}{\usecounter{alphactr}}}{\end{list}}%\alphalist
\def\ul{\underline}
\def\ol{\overline}
\documentclass[12pt]{article}\usepackage{amssymb}
\newtheorem{theorem}{Theorem}[section]
 \newtheorem{lemma}[theorem]{Lemma}%\lemma
 \newtheorem{prop}[theorem]{Proposition}%\prop
\def\ifundefined#1{\expandafter\ifx\csname#1\endcsname\relax}

\begin{document}
\vskip 0.15in
\def\sminus{2 \left[ 1 + {s-1 \choose 3} \right]}

\begin{center}
\bf\Large On sextic surfaces having only nodes\\
(preliminary report)
\end{center}

\vskip 0.15in
\nsfatnebraska

\vspace{0.25in}

\block{Introduction}

Over $\C$, let $S \IN \PP^3$ be a surface having only nodes as singularities.
Let \mp[[ \pi || \tS || S ]] be a minimal resolution of singularities.
A set $N$ of nodes on $S$ is {\it even\/} if there exists a divisor $Q$ on
$\tS$ such that $2Q \sim \pi^{-1}(N)$.

Suppose that $S$ has degree $6$.  It is known (Basset) that $S$ cannot have
$67$ or more nodes.  It is also known (Barth) that $S$ can have $65$ nodes.
It is not known if $S$ can have $66$ nodes.  Likewise, it is not known exactly
what sizes can occur for an even set of nodes on $S$.

We show that an nonempty even set of nodes on $S$ must have size
$24$, $32$, $40$, $56$, or $64$.  We do not know if the sizes $56$ and $64$
can occur.  We show that if $S$ has $66$ nodes, then it must have an even
set of $64$ nodes%
,\footnote{J.\ Wahl informed me that he has independently obtained
these results, but I have not seen his methods.}
and it cannot have an even set of $56$ nodes.  {\bf\it Thus if one could
rule out the case of a 64 node even set, it would follow that S cannot
have 66 nodes.}

The existence or nonexistence of large even node sets is related to the
following vanishing problem.
Let $S \IN \PP^3$ be a normal surface of degree $s$.  Let $D$ be a Weil
divisor on $S$ such that $D \sim_\Q rH$, for some $r \in \Q$.  Under what
circumstances do we have $H^1(\O_S(D)) = 0$?  For instance, this holds when
$r < 0$.  For $s=4$ and $r=0$, $H^1$ can be nonzero.  For $s=6$ and $r=0$, if
a $56$ or $64$ node even set exists, then $H^1$ can be nonzero.  The
vanishing of $H^1$ is also connected to linear normality, quadric normality,
etc.\ of \STCIs\ in $\PP^3$.  See
(\Lcitemark 5\Rcitemark \ 11.8--11.11).  There is no further
discussion of these issues in this paper.

\block{Basic tools and notation}

Since some of the calculations in this report apply as well to surfaces of
arbitrary degree, we will for a while allow $S$ to be any surface of
degree $s$ in $\PP^3$ having only nodes as singularities.

Many ideas used here are taken from Beauville\Lspace \Lcitemark 3\Rcitemark
\Rspace{},
who analyzes quintic surfaces having only nodes.  The following result
is well-known:

\begin{lemma}\label{torsion-free}
$\Pic(\tS)$ is torsion-free.
\end{lemma}

\begin{proof}
Let $\shL$ be a line bundle on $\tS$ such that $\shL^{\o* n} \iso \O_\tS$ for
some $n \in \N$.  Let $S_0$ be the disjoint union of the spectra of the local
rings of the singular points of $S$, and let $\tS_0 = \pi^{-1}(S_0)$.  Since
$\shL$ restricted to each exceptional curve is trivial, it follows
(see\Lspace \Lcitemark 7\Rcitemark \Rspace{}\ \S10, p.\ 157) that
$\shL|_{\tS_0}$ is trivial.  Therefore there exist effective divisors
$D_1, D_2$ on $\tS$ which do not meet the exceptional curves and are such
that $\shL \iso \O_\tS(D_1-D_2)$.  Hence $\shL$ is the pullback of a line
bundle from $S$.  But $\Pic(S)$ is torsion-free and $\Pic(\pi)$ is
injective, so $\shL \iso \O_\tS$.  \qed
\end{proof}

The even sets of nodes on $S$ comprise the codewords of a binary linear
code $C$.  That is, if $C$ is the set of all even sets of nodes, then $C$
may be regarded as a sub-vector-space of $\F_2^n$, where $n$ is the number of
nodes on $S$.  We use without much explanation some standard tools from
coding theory to investigate $C$; the standard reference is
\Lcitemark 6\Rcitemark \Rspace{}.

For $w \in C$, let $\abs{w}$ denote the size of the node set, which in
coding theory would be called the {\it weight\/} of $w$.  It is the number
of $1$'s which appear in the expression for $w$ as an element of $\F_2^n$.
For $n \in \N$, let $C_n$ denote the number of words of weight $n$ in $C$.
Let $C^\perp$ denote the dual code to $C$.

If $J$ is a set of positive integers, by a $[n,k,J]$ code we mean a code
with lives in $\F_2^n$, has dimension $k$, and has nonzero weights in the
set $J$.  We also underscore weights to indicate that they must occur.
Thus, for example, a $[66,13,\setof{24,32,40,\ul{64}}]$ code is a code which
lives in $\F_2^{66}$, has dimension $13$, has nonzero weights in the set
\setof{24,32,40,64}, and has a word of weight $64$.

For $w \in C$, we let $Q_w$ denote a divisor on $\tS$ with the property
that $2Q_w$ is rationally equivalent to the sum of the exceptional curves
corresponding to the points in $w$.  By \pref{torsion-free}, $Q_w$ is
determined up to rational equivalence.

For $w \in C - \setof{0}$, let \mp[[ \tau_w || X_w || \tS ]] be the
double cover branched along the union of the exceptional curves corresponding
to the singular points in $w$, and determined by $Q_w$.
(See\Lspace \Lcitemark 2\Rcitemark \Rspace{}\ I\ \S17 for background on this.)
Since
anything in $\Ker[\Pic(\tau_w)]$ must be $2$-torsion
(\Lcitemark 2\Rcitemark \ I\ 16.2), it follows from \pref{torsion-free} that
$\Pic(\tau_w)$ is injective.

For any divisor $D$ on $\tS$, we can push $D$ down to a Weil divisor $\ol{D}$
on $S$; the exceptional components of $D$ (if any) are ignored.

\block{Results on linear systems and curves}

If $\shL$ is a globally generated line bundle on a smooth projective variety
$T$, we consider the complete linear system corresponding to
$\shL$, and let $\varphi_\shL$ denote the induced morphism from $T$ to its
image, which is a projective variety.  The following result is no
doubt well-known; we include it for lack of a reference.

\begin{lemma}\label{birational}
Let \mp[[ f || W || T ]] be a surjective morphism of smooth projective
varieties of dimension $\geq 2$, which is generically two-to-one, and such that
$\Pic(f)$ is injective.  Let $\shL$ be a globally generated line bundle
on $T$.  Assume that $h^0(f^*\shL) > h^0(\shL)$.  Assume that
$\varphi_\shL$ is birational.  Then $\varphi_{f^*\shL}$ is birational.
\end{lemma}

\begin{proof}
Let $D$ be an effective divisor in the complete linear system associated to
$\shL$.
Let $Q$ be a general member of $\abs{f^*D}$.  By Bertini's theorem, $Q$ is
smooth and connected.  Since $\dim \abs{f^*D} > \dim \abs{D}$, we see that
$Q \notin f^*\abs{D}$.  Let $f_*(Q)$ denote $f(Q)$, thought of as a reduced
effective divisor.  The divisor $f^*f_*(Q) - Q$ is effective.  If
$f^*f_*(Q) = Q$, then $f^*(\shL) \iso f^*\O_T(f_* Q)$, so
$\shL \iso \O_T(f_*Q)$ [since $\Pic(f)$ is injective], and hence
$Q \in f^*\abs{D}$: contradiction.  Hence $f^*f_*(Q) - Q \not= 0$.  Since $Q$
is general it follows that $Q \subsetneq f^{-1}(f(Q))$ as sets.  Therefore if
$x \in Q$ is a general point,
$\setof{x} \subsetneq f^{-1}(f(x))$.  Since $Q$ passes
through $x$ and there is exactly one other point $y \in f^{-1}(f(x))$, which
$Q$ does not pass through, we see that
the restriction of $\varphi_{f^*\shL}$ to $f^{-1}(f(x))$ is
injective.  Since $x$ ranges over general points on general members of
$\abs{f^*D}$, in fact $x$ ranges over a nonempty open subset of $W$.  Hence
there is a nonempty open subset of $W$ on which $\varphi_{f^*\shL}$ is
injective.  \qed
\end{proof}

We state Castelnuovo's bound for the genus of a curve, in a form which
follows immediately from (\Lcitemark 1\Rcitemark \ p.\ 116).

\begin{theorem}\label{castelnuovo}
Let $M$ be a smooth curve of genus $g$, that admits a birational mapping
onto a nondegenerate curve of degree $d$ in $\PP^r$, $r \geq 2$.  Let
$x = \floor{{d-1 \over r-1}}$, $y = {d-1 \over r-1} - x$.
Then:
$${g \over r-1} \leq {x \choose 2} + xy.$$
\end{theorem}

\block{Results valid for any degree {\it s}}

First we make an elementary calculation, and then we give a very useful
proposition, which is probably known in some sense.

\begin{lemma}\label{chi}
Let $w \in C - \setof{0}$, and write $X = X_w$ and $Q = Q_w$ for simplicity.
Let $p = \abs{w}$.  Let $n \in \Z$.  Then
$$\chi(\O_X(n))\ =\ n(n+4-s)s + \sminus - p/4.$$
\end{lemma}

\begin{proof}
\begin{eqnarray*}
\chi(\O_X(n)) & = & \chi(\O_\tS(n)) + \chi(\O_\tS(-Q)(n))\\
& = & \half(nH) \cdot (nH-K_S) + \half(-Q + nH) \cdot (-Q + nH - K_S)
      + 2\chi(\O_\tS)\\
& = & (nH)\cdot(n+4-s)H + 2\chi(\O_\tS) + \half Q^2\\
& = & n(n+4-s)s + \sminus + {1\over8}(\svec E1p)^2.\kern1.15cm\square
\end{eqnarray*}
\end{proof}

Let $w \in C - \setof{0}$, and simplify notation by writing
$Q = Q_w$, $X = X_w$, $\tau = \tau_w$.  Let $E$ be the union of the
exceptional curves corresponding to the singular points in $w$.
For any $n \in \Z$, consider the exact sequence
\sesdot{\O_\tS(nH-Q)}{\O_\tS(nH+Q)}{\O_E(nH+Q)%
}Since the restriction of $nH+Q$ to any component of $E$ has degree $-1$,
$h^j(\O_E(nH+Q)) = 0$ for all $j$.  Hence
$h^k(\O_\tS(nH-Q)) = h^k(\O_\tS(nH+Q))$
for all $k$.  This fact will be used a number of times.

\begin{prop}\label{mmmmm}
Assume that $s \geq 5$.  Assume that $h^0(\O_X(1)) \geq 5$.  Then
$h^0(\O_X(1)) = 5$ and $\abs{w} \leq s(s-1)$.
\end{prop}

\begin{proof}
Let $p = \abs{w}$.  Since $h^0(\O_{\tS}(1)) = h^0(\O_S(1)) = 4$, it follows by
\pref{birational} that $\varphi_{\O_X(1)}$ is birational.  Let $H$ be a
general member of the complete linear system on $\tS$ corresponding to
$\O_{\tS}(1)$.  Let $M = \tau^{-1}(H)$.  By Bertini's theorem, $M$ is a
smooth connected curve.
Let $\O_M(1)$ denote the pullback of $\O_X(1)$ to $M$.
Since $\varphi_{\O_X(1)}$ is birational, $\varphi_{\O_M(1)}$ is birational.

Let $r = h^0(\O_M(1)) - 1$.  The map \mapx[[ M || H ]] is an \'etale double
cover.  Since $H$ has genus $(s-1)(s-2)/2$, $M$ has genus $(s-1)(s-2)-1$.
Since the image of $H$ in $\PP^3$ has degree $s$, it follows that
$\varphi_{\O_M(1)}(M)$ is a nondegenerate curve of degree $2s$ in $\PP^r$.

We have $r \geq 3$; we will show that in fact $r = 3$.  Apply
\pref{castelnuovo}.  We assume that $s \geq 8$, leaving the cases
$s = 5,6,7$ to the reader.  We have
$x = \floor{{2s-1 \over r-1}}$.  Since $y \leq (r-2)/(r-1)$, we have:
\begin{eqnarray*}
{g \over r-1} & = & {(s-1)(s-2)-1 \over r-1}\ \leq\ {x(x-1) \over 2} +
  \left( {r-2 \over r-1} \right) x\\
& \leq & {(2s-1)(2s-r) \over 2(r-1)^2} + {(2s-1)(r-2) \over (r-1)^2},
\end{eqnarray*}
which implies that
$$2(r-1)(s^2 - 3s + 1)\ \leq\ (2s-1)(2s+r-4).$%
$In particular, if $r \geq 4$, then the above inequality holds even when
$r = 4$ is formally substituted in, and it follows that $s^2 - 8s + 3 \leq 0$.
Hence $s < 8$: contradiction.  Hence $r = 3$.

{}From the exact sequence
\ses{\O_X}{\O_X(1)}{\O_M(1)%
}we get $h^0(\O_X(1)) \leq 5$, so in fact $h^0(\O_X(1)) = 5$.  Since
$h^0(\O_X(1)) = h^0(\O_{\tS}(H)) + h^0(\O_{\tS}(H-Q))$, we get
$h^0(\O_{\tS}(H-Q)) = 1$.

Let $I \in \abs{H-Q}$.  Since $I \cdot H = s$, $\pi(\red{I})$ is a curve
(possibly reducible) of degree $\leq s$ on $S$, which passes through at least
$p$ of its nodes.  It follows that $p \leq s(s-1)$.  \qed
\end{proof}

\block{Words of weight < 24 in {\it C}}

{\bf From now on we will assume that {\it s} = 6.}  We will want to get as
much information as possible about the code $C$.  First, an argument of Reid
(\Lcitemark 4\Rcitemark \ 2.11) shows
that all the weights of words in $C$ are divisible by $8$.

Suppose that $w \in C$ is a nonzero word of weight $< 24$.  Then
$\abs{w} \leq 16$.  Let $X = X_w$.  By \pref{chi}, $\chi(\O_X(1)) \geq 12$.
By Serre duality (roughly), $h^0(\O_X(1)) = h^2(\O_X(1))$, so
$h^0(\O_X(1)) \geq 6$, contradicting \pref{mmmmm}.  Hence
$C$ has no nonzero words of weight $< 24$.

\block{Properties of words of weight 24 in {\it C}}\label{weight-24-section}

Let $w \in C_{24}$, and let $X = X_w$.
By \pref{chi}, $\chi(\O_X(1)) = 10$, and since
$h^0(\O_X(1)) = h^2(\O_X(1))$ (as above), we have $h^0(\O_X(1)) \geq 5$.
By \pref{mmmmm}, $h^0(\O_X(1)) = 5$.  Hence $h^0(\O_\tS(H-Q_w)) = 1$.
In other words $\abs{H-Q_w}$ has a unique element, which we henceforth
denote by $D_w$.

We will show now that $\ol{D_w}$ is reduced.  Since $\ol{D_w}$ passes through
$\geq 24$ singular points of $S$, $\deg(\red{(\ol{D_w})}) \geq 5$.  Suppose
that $\ol{D_w}$ is nonreduced.  Then $\deg(\red{(\ol{D_w})}) = 5$ so $\ol{D_w}$
contains a double line.  Let $M$ be the remainder of $\ol{D_w}$, sans the
doubled line.  Since $\ol{D_w}$ must pass with
{\it odd multiplicity\/} through exactly $24$ singular points of $S$, we
conclude that $M$ must pass with odd multiplicity through exactly $24$
singular points of $S$.  Since $\deg(M) = 4$, this is impossible.  Hence
$\ol{D_w}$ is reduced.

Now there is a ``quadric surface'' $V_w \IN \PP^3$ such that
$2\ol{D_w} = V_w \cap S$ as Weil divisors on $S$.  There are four
possibilities for $V_w$.   The worst case, that $V_w$ might be a doubled plane,
can be ruled out, since then $\ol{D_w} = \red{(V_w)} \cap S$ as
Weil divisors on $S$, so $\ol{D_w}$ is Cartier on $S$, which yields a
contradiction.  Unfortunately, we do not know if $V_w$ can be the union of
two distinct planes, a possibility which substantially complicates the
proofs.  For starters, we must analyze this case.  We could rule it out
if we could show that a cubic curve on $S$ cannot pass through $15$ singular
points of $S$.

\begin{prop}\label{two-planes}
Let $w \in C_{24}$.  Suppose that $V_w$ is the union of two distinct planes
$H_1$, $H_2$.  For each $i$, let $\ol{D_i} = \red{(H_i \cap S)}$.  Then:
\begin{alphalist}
\item for each $i$, $\ol{D_i}$ passes through exactly $15$ singular points of
      $S$, and it passes with multiplicity one through each of these;
\item of the $15$ singular points, exactly $12$ are in $w$, and the remaining
      $3$ lie on the line $H_1 \cap H_2$.  Both $\ol{D_1}$ and $\ol{D_2}$
      pass with multiplicity one through these $3$ nodes.
\end{alphalist}
\end{prop}

\begin{proof}
Let $D_i \IN \tS$ be the strict transform of $\ol{D_i}$, for each $i$.  Let
$L = H_1 \cap H_2$.  Because $2|\deg(S)$ and $\ol{D_w}$ is reduced, we cannot
have $L \IN S$.  Hence $D_1$ meets $D_2$ properly.  Note that $\deg(D_i) = 3$.

Let $Q_i = \half \sum_E (D_i \cdot E) E$, as $E$ ranges over the exceptional
curves of $\tS$.  Then $Q$ is a $\Q$-divisor, but (caution) it is not a
divisor.  We have $D_i \sim_\Q \half H - Q_i$.

Suppose that $\ol{D_1}$ has only nodes as singularities.  Then
$\chi(\O_{D_1})$ is equal to the number (say $n$) of such singular points
which are resolved by $\pi$.  Then $D_1^2 = -2n-6$, from which it
follows that $Q_1$ has the form
$$\half\left(\svec E1{15} + \sum_{i=1}^n 2E_i'\right),$%
$where $\vec E1{15}, E_1',\ldots,E_n'$
are distinct exceptional curves.  But $\ol{D_1}$
has degree $3$, and so can pass through at most $3(\deg(S)-1) = 15$ singular
points of $S$.  Hence $n = 0$, so (a) holds, for $D_1$.

In the three cases where $\ol{D_1}$ has a singularity other than a node,
arguments similar to those of the last paragraph still yield the conclusion
that (a) holds, for $D_1$.

Of course, the results of the last two paragraphs apply equally to $\ol{D_2}$.
Since $\ol{D_w}$ cannot pass with multiplicity one through a node outside $w$,
it follows that the number of nodes of $w$ (say $r$) which $\ol{D_i}$ passes
through is independent of $i$.

Since $\ol{D_1}$ and $\ol{D_2}$ must together pass with multiplicity one
through all $24$ of the singular points in $w$, we must
have $r = 12$.  Statement (b) follows.  \qed
\end{proof}

\begin{lemma}\label{line-through-four}
Let $L \IN S$ be a line.  Let $w \in C$.  Then $L$ can pass through
at most $4$ singular points of $S$ which are in $w$.
\end{lemma}

\begin{proof}
Let $\tL \IN \tS$ be the strict transform of $L$.
Certainly $L$ can pass through at most $5$ singular points of $S$.
But $\tL \cdot Q_w \in \Z$, so the lemma follows.  \qed
\end{proof}

\begin{prop}\label{empty-intersection}
Let $w_1, w_2 \in C$, and assume that
$\abs{w_1} = \abs{w_2} = \abs{w_1+w_2} = 24$.  Then either (a)
$D_{w_1} \cap D_{w_2} = \emptyset$, or else (b) we can write
$V_{w_1} = H_1 \cup H_2$ and $V_{w_2} = H_2 \cup H_3$, where $H_1, H_2, H_3$
are planes.
\end{prop}

\begin{proof}
First note that $(H-Q_{w_1}) \cdot (H-Q_{w_2}) = 6 - \half(12) = 0$, so for
case (a) it suffices to show that $D_{w_1}$ and $D_{w_2}$ share no component.
Let $M$ be the shared part of $D_{w_1}$ and $D_{w_2}$.

We may assume that $V_{w_1}$ meets $V_{w_2}$ properly, since otherwise
conclusion (b) holds.  Then $\deg(V_{w_1} \cap V_{w_2}) = 4$.  We have
$\oM \IN V_{w_1} \cap V_{w_2}$, and since $V_{w_i}$ is tangent to $S$
along $\oM$ for each $i$, in fact $V_{w_1}$ is tangent to $V_{w_2}$ along
$\oM$.  Hence $2\oM \IN V_{w_1} \cap V_{w_2}$.  Hence $\deg(M) \leq 2$.
If $\deg(M) = 2$, then $V_{w_1} \cap V_{w_2} = 2\oM$.
Hence $\oM$ passes through
all $12$ singular points in $w_1 \cap w_2$.  But $\deg(\oM) = 2$, so $\oM$
can pass through at most $2[\deg(S)-1] = 10$ singular points of $S$:
contradiction.  Hence $\deg(M) \leq 1$.  We have
$$0\ \leq\ (H-Q_{w_1}-M) \cdot (H-Q_{w_2}-M) = M(Q_{w_1}+Q_{w_2}-2H)+M^2.
  \eqno(*)$%
$We will complete the proof by showing that if $M \not= 0$, then the \RHS\ of
$(*)$ is negative.  Write $M = R+E$, where $R$, $E$ are effective, $E$ is
exceptional, and $R$ has no exceptional components.  Expanding out the
\RHS\ of $(*)$ yields
$$(R+E)\cdot(Q_{w_1}+Q_{w_2}) - 2\deg(R) + R^2 + 2(E \cdot R) + E^2.$%
$If $R=0$, then $E=0$, or else this quantity is negative.  Hence $\deg(R) = 1$.

By the adjunction formula, $R^2 = -4$.  Therefore
to get a contradiction, it is enough to show that
$$(R+E) \cdot (Q_{w_1} + Q_{w_2}) + 2(E \cdot R) + E^2 < 6.$%
$By \pref{line-through-four}, $R \cdot (Q_{w_1} + Q_{w_2}) \leq 4$, so it is
enough to show that
$$E \cdot (Q_{w_1} + Q_{w_2}) + 2(E \cdot R) + E^2 < 2.$%
$Thus it is sufficient to show that if $F$ is an exceptional curve, then for
all $n \in \N$,
$$nF \cdot (Q_{w_1} +Q_{w_2}) + 2(nF \cdot R) + (nF)^2 \leq 0.$%
$Clearly it is enough to do the case $n=1$, so we need
$$F \cdot(Q_{w_1} + Q_{w_2}) + 2(F \cdot R) \leq 2.$%
$Since $F \cdot R \leq 1$, this is clear.  \qed
\end{proof}

\block{Words of weight 48 in {\it C}}

We will show that $C$ has no words of weight $48$.  So let
$v \in C_{48}$, heading for a contradiction.  Write $Q$ for $Q_v$, $X$ for
$X_v$.  The first step is to show that
$h^0(\O_\tS(2H-Q)) = 1$.  Let $N \IN \F_2^{48}$ be $\setof{w \in C: w \IN v}$.

We have $h^0(\O_\tS(H-Q)) = 0$, since otherwise there is a curve of degree
$6$ on $S$ which passes through $\geq 48$ singular points of $S$.  Hence
$h^0(\O_X(1)) = 4$.  By \pref{chi} we have $\chi(\O_X) = 10$, so
$h^2(\O_X) - h^1(\O_X) = 9$.  But $h^2(\O_X) \geq h^2(\O_\tS) = 10$, so
$h^1(\O_X) \geq 1$.  Hence $\dim_{\F_2}(\Pic(X)_2) \geq 2$.  Hence by
(\Lcitemark 3\Rcitemark \ Lemma 2), $\dim(N) \geq 3$.  Since we know
that $w \in N-\setof{0}\ \Longrightarrow\ \abs{w} \geq 24$, we in fact have
$\setof{\abs{w}: w \in N} = \setof{0,24,48}$.  Choose $w_1, w_2 \in N$ with
$\abs{w_1} = \abs{w_2} = 24$ and $w_1 \cap w_2 = \emptyset$.  By the first
paragraph of \S\ref{weight-24-section},
$h^0(\O_\tS(H-Q_{w_i})) \not= 0$.  Since $v = w_1 + w_2$,
$h^0(\O_\tS(2H-Q)) \not= 0$.  By Serre duality,
$h^2(\O_\tS(Q)) \not= 0$.  Hence $h^2(\O_\tS(-Q)) \not= 0$.  Then
$h^2(\O_X) = h^2(\O_\tS) + h^2(\O_\tS(-Q)) \geq 11$, so $h^1(\O_X) \geq 2$.
Hence $\dim(N) \geq 5$.

Suppose that $h^0(\O_\tS(2H-Q)) > 1$.  Then (arguing as above)
$h^1(\O_X) \geq 3$, so $\dim(N) \geq 7$.  Hence there exists a
$[48,7,\setof{24,48}]$ code, which contradicts the Griesmer bound for
codes.  (One also gets a contradiction by the linear programming method.)
Hence $h^0(\O_\tS(2H-Q)) = 1$.

Let $w \in N_{24}$.  Since $h^0(\O_\tS(2H-Q)) = 1$, it follows that
$D_w + D_{v+w}$ is independent of $w$.
Since $\dim(N) \geq 3$, we can find $w_1, w_2 \in N$ with
$\abs{w_1} = \abs{w_2} = \abs{w_1+w_2} = 24$.  We will show that
$V_{w_1}$ does not meet $V_{w_2}$ properly.  Suppose otherwise.
We have $D_{w_1} + D_{v+w_1} = D_{w_2} + D_{v+w_2}$.
By \pref{empty-intersection}, $D_{w_1} \cap D_{w_2} = \emptyset$.
Hence $D_{w_1} \IN D_{v+w_2}$, so $\ol{D_{w_1}} = \ol{D_{v+w_2}}$, so
$w_1 = v+w_2$: contradiction.  Hence $V_{w_1}$ does not meet
$V_{w_2}$ properly, i.e.\ each consists of two planes, with one in common.
The same holds for the pair $(V_{w_1}, V_{v+w_2})$, the pair
$(V_{v+w_1}, V_{w_2})$, and the pair $(V_{v+w_1}, V_{v+w_2})$.

We consider now the configuration $\cal C$ of planes which occur in $V_{w_1}$,
$V_{w_2}$, $V_{v+w_1}$, and $V_{v+w_2}$.  In total, there are four planes.

First suppose that no three of these planes share a common line.  Then
we may visualize the configuration by means of a tetrahedron $T$.  (The faces
correspond to planes.)  For any pair of faces, the intersection of the
corresponding two planes with $S$ gives rise to a word $w$ of weight $24$,
and so \pref{two-planes} applies.  Therefore every edge of $T$ has exactly $3$
singular points of $S$ on it, and every face of $T$ has exactly $15$ singular
points of $S$ on it.  Here is the generic picture, showing only the singular
points which appear on the edges:

\vspace{0.1in}

\par\noindent{\bf[There was a postscript picture of a tetrahedron here, with
three little spheres along each edge.]}

\vspace{0.1in}

\par\noindent
The picture is generic because some of the edge singular points could be
at the vertices.  Restricting to the generic case for the moment, let us
compute the total number of singular points of $S$ which lie on $T$.  Since
each face contains $15$, the interior of each face has $6$.  Hence
$$\abs{\Sing(S) \cap T}\ =\ (4 \cdot 6) + (6 \cdot 3)\ =\ 42.$%
$But $T$ must go through at least $\abs{v} = 48$ singular points of $S$:
contradiction.  In the nongeneric cases, one similarly gets a
contradiction.  (In the least generic case, where all $4$ vertices of $T$ are
in $\Sing(S)$, one gets $\abs{\Sing(S) \cap T} = 46$, which still gives a
contradiction.)
Hence at least $3$ of the planes in $\cal C$ must share
a common line.

Suppose now that the four planes in $\shC$ do not all share a common line.
Taking a cross section by a suitable plane, we may represent each plane by
a line segment and each intersection of two planes by a vertex.  (The single
point where all planes meet is not shown.)  In the generic case, we can
label each vertex in the picture with a $3$.

\widepost{85}{25}{
newpath 0 -15 moveto -50 50 lineto stroke
newpath 0 -15 moveto 50 50 lineto stroke
newpath 0 -15 moveto 0 50 lineto stroke
newpath -50 50 moveto 50 50 lineto stroke
0 -15 PointAt -50 50 PointAt 50 50 PointAt 0 50 PointAt
{6 270 (3) AnglePrint} 0 -15 xyput
{6 90 (3) AnglePrint} 0 50 xyput
{6 135 (3) AnglePrint} -50 50 xyput
{6 45 (3) AnglePrint} 50 50 xyput}

\par\noindent One finds that the four planes together pass through at
most $45$ singular points of $S$: contradiction.  In the nongeneric case,
one finds that $\shC$ passes through $47$ points of $\Sing(S)$, still shy
of $48$.

Now we know that the four planes in $\shC$ must contain a common line.
It follows that there exists a line $L$ such that for any
$w \in N_{24}$, $L$ is contained in both planes of $V_w$.
Let $\shC^+$ be the configuration
of all planes which appear in $V_w$ for some $w \in N_{24}$.
Since $\dim(N) \geq 5$, $\abs{N_{24}} \geq 30$.  Since
${8 \choose 2} < 30$, $\shC^+$ must have at least $9$ planes in it.  Hence
$\shC^+$ must pass through at least $(9 \cdot 12) + 3 = 111$ singular points
of $S$, which is absurd.  We conclude that $C$ has no word of weight $48$.

\block{Apply coding theory}

Suppose now that $S$ has $66$ nodes.  Since $H^2(S,\Z) \iso \Z^{106}$, it
follows\Lspace \Lcitemark 3\Rcitemark \Rspace{} that $\dim(C) \geq 66 - 106/2 =
13$.
We know that the nonzero weights appearing in $C$ are all in the set
\setof{24,32,40,56,64}.  Now we
show that weight $64$ must occur and that weight $56$ cannot occur.

\begin{theorem}\label{code-doesnt-exist}
There is no $[66,13,\setof{24,32,40,56}]$ code.
\end{theorem}

\begin{proof}
Let $R = \setof{24,32,40,56}$ and $T = \setof{4,8,12,16,20}$.  By the linear
programming method (LP), there is no $[21,10,T]$ code.  (See the appendix
for some {\sc Maple} code to do calculations like this.)  Therefore by taking
the residual code \WRT\ a codeword of weight $40$ (hereafter referred to as
res!40), one sees that there is no $[61,11,\setof{24,32,\ul{40},56}]$ code.

Let $C$ be a $[58,11,\setof{24,32,56}]$ code.  By LP, $C_{56} = 1$.  By LP,
$C^\perp_1 = 2$.  By LP, $C^\perp_2 = 6.5$: contradiction.  Hence there is no
such code $C$.  Now let $C$ be a $[59,11,\setof{24,32,56}]$ code.  Then
$C^\perp_1 = 0$.  Then LP yields a contradiction, so there is no such code $C$.
Similarly, one sees that there is no $[60,11,\setof{24,32,56}]$ code, and
finally in the same way that $(*)$ there is no $[61,11,\setof{24,32,56}]$ code.
By the preceeding paragraph, we conclude $(\dag)$ that there is no $[61,11,R]$
code.

By LP, there is no $[22,11,T]$ code.  Therefore by res!40, there is no
$[62,12,\setof{24,32,\ul{40},56}]$ code.  By LP, there is no
$[62,12,\setof{24,32,56}]$ code $C$ with $C^\perp_1 = 0$, and therefore by
$(*)$ of the preceeding paragraph, there is no $[62,12,\setof{24,32,56}]$ code
at all.  Hence there is no $[62,12,R]$ code.

{}From $(\dag)$ we see that a $[63,13,R]$ code $C$ must have
$C^\perp_2 = C^\perp_3 = 0$.  By LP, one concludes that there is no
$[63,13,R]$ code.

{}From the preceeding three paragraphs we see that a $[64,13,R]$ code $C$ must
have $C^\perp_1 = C^\perp_2 = C^\perp_3 = 0$.  By LP, one concludes that
a $[64,13,R]$ code $C$ must have $C_{56} = 2.5$: contradiction.  Hence there
is no such code $C$.

Let $C$ be a $[65,13,R]$ code.  Then $C^\perp_1 = 0$.  By LP,
$C^\perp_2 \geq 5$.  Consider two distinct words of weight two in $C^\perp$.
If they are disjoint [resp.\ not disjoint], then there is a
$[61,11,R]$ code [resp.\ $[62,12,R]$ code], in either case contradicting
earlier results.  Hence there is no $[65,13,R]$ code.

Let $C$ be a $[66,13,R]$ code.  (This will lead to a contradiction.)  Then
$C^\perp_1 = 0$.  Form a graph $G$ whose vertices correspond to the coordinate
positions covered by the words of weight $2$ in $C^\perp$, and such that the
edges of $G$ correspond to the words of weight $2$ in $C^\perp$.  Then $G$ is a
disjoint union of complete graphs.  By LP, $C^\perp_2 \geq 7$, so $G$ has at
least $7$ edges.  Now $G$ has no $K_4$ component, because then one would have a
$[62,12,R]$ code.  Similar arguments show that $G$ cannot have a $K_n$ (for
$n \geq 5$), it cannot have $2$ $K_3$'s, and it cannot have a $K_3$ and a
$K_2$.  Therefore, since $G$ has at least $7$ edges, it must consist entirely
of $K_2$'s.  Otherwise said, the words of weight $2$ in $C^\perp$ are pairwise
disjoint.

If $C$ has a word $v$ of weight $56$, it occupies all but $10$ bits of $C$, and
so at most $5$ of the words of weight $2$ in $C^\perp$ are disjoint from $v$.
Thus we can choose some $w \in C^\perp_2$ with the property that if such
a word $v$ exists, then $w \IN v$.

Let $D$ be the subcode of $C$ consisting of words disjoint from $w$.  Then $D$
is a $[64,12,\setof{24,32,40}]$ code and $D^\perp$ has at least $7-1=6$ words
of weight $2$.  Adjoin the unique word of weight $64$ to $D$,
yielding a $[64,13,\setof{24,32,40,64}]$ code $E$ with $E^\perp_2 \not= 0$.
Hence there exists a $[62,12,\setof{24,32,40}]$ code: contradiction.  \qed
\end{proof}

We conclude that $C$ must be a
$[66,13,\setof{24,32,40,\ul{64}}]$ code.  Unfortunately, such a code exists.

\vspace{0.1in}

\par\noindent{\bf Appendix}

\vspace{0.1in}

Here we give some {\sc Maple} code which can be used to carry out the
linear programming calculations which appear in the proof of
\pref{code-doesnt-exist}.

The first four lines can be adjusted to check
for particular types of codes.  As written, the program will analyze
$[64,13,\setof{24,32,40,56}]$ codes $C$ with
$C^\perp_1 = C^\perp_2 = C^\perp_3 = 0$.  If you
run it, you will find that such a code must have $C_{56} = 2.5$ (which is
absurd), as claimed in the proof of \pref{code-doesnt-exist}.

\vspace{0.2in}

\begin{verbatim}
n := 64;
k := 13;
w := {24,32,40,56};
extraconstraints := [mu1=0, mu2=0, mu3=0];

K := (x,i) -> sum( (-1)^j * binomial(i,j) * binomial(n-i,x-j), j=0..x ):

for m from 0 to n/2 do
     mu.m := (K(m,0) + convert(map(i -> 'a.i' * K(m,i), w), `+`)) / 2^k:
od:

objective := 1 + convert(map(z -> a.z, w), `+`);

mus := {seq( mu.m >= 0, m=0..n/2 ), objective = 2^k, op(extraconstraints)}:

for i from 1 to 2 do
     M := simplex[minimize](mu.i,mus,NONNEGATIVE);
     `minimum value of mu`.i, evalf(subs(M,mu.i));
     M := simplex[maximize](mu.i,mus,NONNEGATIVE);
     `maximum value of mu`.i, evalf(subs(M,mu.i));
od;

for q in w do
     M := simplex[minimize](a.q,mus,NONNEGATIVE);
     q, `min`, evalf(subs(M,a.q));
     M := simplex[maximize](a.q,mus,NONNEGATIVE);
     q, `max`, evalf(subs(M,a.q));
od;

quit():
\end{verbatim}
\message {DOCUMENT TEXT}
\section*{References}
\addcontentsline{toc}{section}{References}
\ \par\noindent\vspace*{-0.25in}
\hfuzz 5pt

\message{REFERENCE LIST}

\bgroup\Resetstrings%
\def\Loccittest{}\def\Abbtest{ }\def\Capssmallcapstest{}\def\Edabbtest{
}\def\Edcapsmallcapstest{}\def\Underlinetest{ }%
\def\NoArev{1}\def\NoErev{0}\def\Acnt{4}\def\Ecnt{0}\def\acnt{0}\def\ecnt{0}%
\def\Ftest{ }\def\Fstr{1}%
\def\Atest{ }\def\Astr{Arbarello\Revcomma E\Initper %
  \Acomma M\Initper  Cornalba%
  \Acomma P\Initper \Initgap A\Initper  Griffiths%
  \Aandd J\Initper  Harris}%
\def\Ttest{ }\def\Tstr{Geometry of Algebraic Curves, Vol.\ I}%
\def\Itest{ }\def\Istr{Spring\-er-Ver\-lag}%
\def\Ctest{ }\def\Cstr{New York}%
\def\Dtest{ }\def\Dstr{1985}%
\Refformat\egroup%

\bgroup\Resetstrings%
\def\Loccittest{}\def\Abbtest{ }\def\Capssmallcapstest{}\def\Edabbtest{
}\def\Edcapsmallcapstest{}\def\Underlinetest{ }%
\def\NoArev{1}\def\NoErev{0}\def\Acnt{3}\def\Ecnt{0}\def\acnt{0}\def\ecnt{0}%
\def\Ftest{ }\def\Fstr{2}%
\def\Ven{Van de Ven}{}%
\def\Atest{ }\def\Astr{Barth\Revcomma W\Initper %
  \Acomma C\Initper  Peters%
  \Aandd A\Initper  \Ven}%
\def\Ttest{ }\def\Tstr{Compact Complex Surfaces}%
\def\Itest{ }\def\Istr{Spring\-er-Ver\-lag}%
\def\Ctest{ }\def\Cstr{New York}%
\def\Dtest{ }\def\Dstr{1984}%
\def\Qtest{ }\def\Qstr{access via "barth peters ven"}%
\def\Xtest{ }\def\Xstr{local-triv: I 10.1 (theorem of Grauert-Fischer)}%
\Refformat\egroup%

\bgroup\Resetstrings%
\def\Loccittest{}\def\Abbtest{ }\def\Capssmallcapstest{}\def\Edabbtest{
}\def\Edcapsmallcapstest{}\def\Underlinetest{ }%
\def\NoArev{1}\def\NoErev{0}\def\Acnt{1}\def\Ecnt{1}\def\acnt{0}\def\ecnt{0}%
\def\Ftest{ }\def\Fstr{3}%
\def\Atest{ }\def\Astr{Beauville\Revcomma A\Initper }%
\def\Ttest{ }\def\Tstr{Sur le nombre maximum de points doubles d'une surface in
$\PP^3$  $(\mu(5) = 31)$}%
\def\Btest{ }\def\Bstr{Algebraic Geometry, Angers, 1979}%
\def\Etest{ }\def\Estr{A\Initper  Beauville}%
\def\Itest{ }\def\Istr{Sijthoff \& Noordhoff}%
\def\Dtest{ }\def\Dstr{1980}%
\def\Ptest{ }\def\Pcnt{ }\def\Pstr{207--215}%
\def\Qtest{ }\def\Qstr{access via "beauville double points"}%
\Refformat\egroup%

\bgroup\Resetstrings%
\def\Loccittest{}\def\Abbtest{ }\def\Capssmallcapstest{}\def\Edabbtest{
}\def\Edcapsmallcapstest{}\def\Underlinetest{ }%
\def\NoArev{1}\def\NoErev{0}\def\Acnt{1}\def\Ecnt{0}\def\acnt{0}\def\ecnt{0}%
\def\Ftest{ }\def\Fstr{4}%
\def\Atest{ }\def\Astr{Catanese\Revcomma F\Initper }%
\def\Ttest{ }\def\Tstr{Babbage's conjecture, contact of surfaces, symmetric
determinantal varieties and applications}%
\def\Jtest{ }\def\Jstr{Invent. Math.}%
\def\Vtest{ }\def\Vstr{63}%
\def\Dtest{ }\def\Dstr{1981}%
\def\Ptest{ }\def\Pcnt{ }\def\Pstr{433--465}%
\def\Qtest{ }\def\Qstr{access via "catanese contact"}%
\Refformat\egroup%

\bgroup\Resetstrings%
\def\Loccittest{}\def\Abbtest{ }\def\Capssmallcapstest{}\def\Edabbtest{
}\def\Edcapsmallcapstest{}\def\Underlinetest{ }%
\def\NoArev{1}\def\NoErev{0}\def\Acnt{1}\def\Ecnt{0}\def\acnt{0}\def\ecnt{0}%
\def\Ftest{ }\def\Fstr{5}%
\def\Atest{ }\def\Astr{Jaffe\Revcomma D\Initper \Initgap B\Initper }%
\def\Ttest{ }\def\Tstr{Applications of iterated curve blowup to set theoretic
complete intersections in $\PP^3$}%
\def\Rtest{ }\def\Rstr{preprint}%
\def\Qtest{ }\def\Qstr{access via "jaffe rational singularities" or preferably
"jaffe applications intersections"}%
\def\Htest{ }\def\Hstr{7}%
\Refformat\egroup%

\bgroup\Resetstrings%
\def\Loccittest{}\def\Abbtest{ }\def\Capssmallcapstest{}\def\Edabbtest{
}\def\Edcapsmallcapstest{}\def\Underlinetest{ }%
\def\NoArev{1}\def\NoErev{0}\def\Acnt{2}\def\Ecnt{0}\def\acnt{0}\def\ecnt{0}%
\def\Ftest{ }\def\Fstr{6}%
\def\Atest{ }\def\Astr{MacWilliams\Revcomma F\Initper \Initgap J\Initper %
  \Aand N\Initper \Initgap J\Initper \Initgap A\Initper  Sloane}%
\def\Ttest{ }\def\Tstr{The Theory of Error-Correcting Codes}%
\def\Itest{ }\def\Istr{North-Holland}%
\def\Ctest{ }\def\Cstr{Amsterdam}%
\def\Dtest{ }\def\Dstr{1977}%
\def\Qtest{ }\def\Qstr{access via "macwilliams sloane"}%
\Refformat\egroup%

\bgroup\Resetstrings%
\def\Loccittest{}\def\Abbtest{ }\def\Capssmallcapstest{}\def\Edabbtest{
}\def\Edcapsmallcapstest{}\def\Underlinetest{ }%
\def\NoArev{1}\def\NoErev{0}\def\Acnt{1}\def\Ecnt{0}\def\acnt{0}\def\ecnt{0}%
\def\Ftest{ }\def\Fstr{7}%
\def\Atest{ }\def\Astr{Pinkham\Revcomma H\Initper }%
\def\Ttest{ }\def\Tstr{Singularites rationnelles de surfaces}%
\def\Btest{ }\def\Bstr{S\'eminaire sur les Singularit\'es des Surfaces}%
\def\Stest{ }\def\Sstr{Lecture Notes in Mathematics}%
\def\Itest{ }\def\Istr{Spring\-er-Ver\-lag}%
\def\Ctest{ }\def\Cstr{New York}%
\def\Vtest{ }\def\Vstr{777}%
\def\Dtest{ }\def\Dstr{1980}%
\def\Ptest{ }\def\Pcnt{ }\def\Pstr{147--178}%
\def\Qtest{ }\def\Qstr{access via "pinkham rational singularities"}%
\Refformat\egroup%

\end{document}